\documentclass{amsproc}

\usepackage{}

\newtheorem{theorem}{Theorem}[section]
\newtheorem{proposition}[theorem]{Proposition}

\theoremstyle{definition}

\newtheorem{examples}[theorem]{Examples}
\newtheorem{hypothesis}[theorem]{Hypothesis}

\theoremstyle{remark}
\newtheorem{remark}[theorem]{Remark}

\newcommand{\HH}{\mathrm{H}}

\newcommand{\cat}{\hat{\mathcal{C}}}
\newcommand{\Aut}{\mathrm{Aut}}
\newcommand{\Def}{\mathrm{Def}}

\numberwithin{equation}{section}

\begin{document}

\title[Universal deformation rings of group representations]{Universal deformation rings of group 
representations, with an application of Brauer's generalized decomposition numbers}

\author{Frauke M. Bleher}
\address{Department of Mathematics\\University of Iowa\\
Iowa City, IA 52242-1419}
\email{frauke-bleher@uiowa.edu}
\thanks{The author was supported in part by  
NSA Grant H98230-11-1-0131.}

\begin{abstract}
We give an introduction to the deformation theory of linear representations of profinite groups
which Mazur initiated in the 1980's. We then consider the case of representations of finite groups. 
We show how Brauer's generalized decomposition numbers can be used in some cases to
explicitly determine universal deformation rings.
\end{abstract}

\maketitle


\section{Introduction}
\label{s:intro}

Let $p$ be a prime number and let $(K,\mathcal{O},k)$ be a $p$-modular system, where 
$\mathcal{O}$ is a complete discrete valuation ring of characteristic 0 with maximal
ideal $\mathfrak{m}_\mathcal{O}$, $K=\mathrm{Frac}(\mathcal{O})$ is its fraction field 
and $k=\mathcal{O}/\mathfrak{m}_\mathcal{O}$ is its residue field of characteristic $p$.
Let $G$ be a group and let $V$ be a $kG$-module of finite $k$-dimension.
It is a classical question to ask whether $V$ can be lifted to an $\mathcal{O}$-free
$\mathcal{O}G$-module. Green showed in \cite{green} that this is possible if
$\mathrm{Ext}^2_{kG}(V,V)=0$. However, this is only a sufficient criterion, as there
are many cases when $\mathrm{Ext}^2_{kG}(V,V)\neq 0$ and $V$ can still be lifted
to an $\mathcal{O}$-free $\mathcal{O}G$-module. Even if we know that there is such a lift of $V$
to $\mathcal{O}$, we may still not be able to characterize all possible lifts to $\mathcal{O}$.
On the other hand, if $V$ cannot be lifted to $\mathcal{O}$, it is natural to ask if $V$ can 
be lifted to other complete local commutative rings with residue field $k$. 
This can be formulated as the following two natural questions:
\begin{itemize}
\item[(i)] How can all possible lifts of $V$ to $\mathcal{O}$ be described?
\item[(ii)] Over which complete local commutative rings with residue field $k$ can $V$ be lifted?
Is there one particular such complete local ring from which all these lifts arise?
\end{itemize}
To answer both of these questions, it is necessary to develop a systematic way
to study isomorphism classes of lifts, also called deformations, of $kG$-modules $V$. 
In the 1980's, Mazur developed a deformation theory of representations of profinite
Galois groups over finite fields
which can be used to give answers to these questions in certain cases.

The goal of this paper is to give an introduction to this deformation theory,
with an emphasis on the particular case of deformations of representations of finite
groups. 

In section \ref{s:maz}, we will first focus on Mazur's  general deformation theory
by studying the deformation functor, universal deformation rings, tangent spaces,
obstructions, and deformation rings in number theory. 
In section \ref{s:udrfinite}, we will then concentrate on deformation rings and
deformations of representations of finite groups. In section
\ref{s:brauer}, we will describe how Brauer's generalized decomposition numbers
can be used to explicitly determine the universal deformation rings of
certain representations belonging to blocks of tame representation type.

We close this introduction by discussing
a few classical examples which show that lifts 
occur naturally both in representation theory and number theory.
 
\begin{examples}
\hspace*{1cm}
\begin{enumerate}
\item[(i)] Let $G$ be a finite group and consider permutation modules for $kG$. 

Recall that
direct summands of permutation modules are also called $p$-permutation modules or
trivial source modules. Scott proved in \cite{scott} that $p$-permutation modules
can be lifted to $\mathcal{O}$. Rickard used this result, for example, in \cite{splendid} to show that
tilting complexes defining splendid equivalences can be lifted from $k$ to $\mathcal{O}$.

Another class of permutation modules is given by endo-permutation modules for
$kG$, i.e. $kG$-modules $V$ for which $\mathrm{End}_k(V)$ is a permutation module.
A special subclass is provided by the endo-trivial modules, for which $\mathrm{End}_k(V)$
is isomorphic as a $kG$-module to a direct sum of the trivial simple $kG$-module and
a projective $kG$-module. Alperin proved in \cite{alperinendotrivial} that if $G$ is a 
$p$-group and $V$ is an endo-trivial $kG$-module, then $V$ can be lifted to an 
endo-trivial $\mathcal{O}G$-module.

\vspace{1ex}

\item[(ii)] Let $k=\mathbb{F}_p=\mathbb{Z}/p$ and let $\mathcal{O}=\mathbb{Z}_p$ be the ring of $p$-adic
integers. Suppose $E$ is an elliptic curve over $\mathbb{Q}$ and define $G=G_\mathbb{Q}=
\mathrm{Gal}(\overline{\mathbb{Q}}/\mathbb{Q})$. Then $G$ acts on the 
torsion points $E[p^n]\cong \mathbb{Z}/p^n\times \mathbb{Z}/p^n$ for all $n\ge 1$. 
In particular, for $n=1$, we obtain a $kG$-module $V$ of $k$-dimension 2.
Since the action of $G$ on $E[p^n]$ commutes with the multiplication by $p$ on 
$E[p^n]$, $G$ acts naturally on the Tate module $T_p(E)=\displaystyle\lim_{\longleftarrow}E[p^n]
\cong \mathbb{Z}_p\times\mathbb{Z}_p$. Hence $T_p(E)$ defines a lift of $V$
over $\mathcal{O}=\mathbb{Z}_p$. 
\end{enumerate}
\end{examples}


\section{Mazur's deformation theory}
\label{s:maz}

In the 1980's, Mazur developed a deformation theory of Galois representations
to systematically study $p$-adic
lifts of representations of profinite Galois  groups \cite{maz1}. His initial motivation 
seems to have come from Hida's work on ordinary $p$-adic modular forms 
(see for example Gouv\^{e}a's summary in \cite[Lecture 2]{gouvea2001}).
Mazur's deformation theory became a fundamental tool in number theory when it
was shown to play a crucial role in the proof of Fermat's Last Theorem, and more 
generally in the proof of the modularity conjecture for elliptic curves over $\mathbb{Q}$ 
(see \cite{wi,wita,bcdt}). For more information and background, we also refer the reader to the survey articles
\cite{gouvea1991,gouvea2001,maz2} on Mazur's deformation theory and the
collection of articles in \cite{CMS17,cornell} on the proof of Fermat's Last Theorem by
Wiles and Taylor.

\subsection{The deformation functor}
\label{ss:deffunctor}

Let $k$ be a perfect field of positive characteristic $p$, and let $W=W(k)$ be the ring of
infinite Witt vectors over $k$. Recall that since $k$ is perfect, $W$ is the unique (up to isomorphism)
complete discrete valuation ring of characteristic zero which is absolutely unramified, 
meaning that the valuation of $p$ is 1 so that $p$ is a generator of the maximal ideal 
$\mathfrak{m}_W$. Let $\cat$ denote the category whose objects are all 
complete local commutative Noetherian rings with residue field $k$
and whose morphisms are local homomorphisms of complete local Noetherian rings
that induce the identity on $k$. Strictly speaking, each object of $\cat$ is a pair
$(R,\pi_R)$ consisting of a complete local commutative Noetherian ring $R$ and a fixed
reduction map $\pi_R:R\to k$ inducing an isomorphism $R/\mathfrak{m}_R\cong k$.
Note that all rings $R$ in $\cat$ have a natural $W$-algebra structure, which means that
the morphisms in $\cat$ can also be viewed as continuous $W$-algebra homomorphisms
inducing the identity on the residue field $k$.

Let $G$ be a profinite group, and let $V$ be a finite dimensional vector space over $k$
with a continuous $k$-linear action of $G$ on $V$, which is given by 
a continuous group homomorphism from $G$ to the discrete group $\Aut_k(V)$.
A \emph{lift} of $V$ over a ring $R$ in $\hat{\mathcal{C}}$ is defined to be a pair $(M,\phi)$
consisting of a finitely generated free $R$-module $M$ on which $G$ acts
continuously together with a $G$-equivariant isomorphism 
$\phi: k \otimes_R M \to V$ of $k$-vector spaces.  
We say two lifts $(M,\phi)$ and $(M',\phi')$ of $V$ over $R$ are \emph{isomorphic} if 
there exists an $R$-linear $G$-equivariant isomorphism $f:M\to M'$ satisfying
$\phi'\circ(k\otimes f)=\phi$. 
We define the set $\mathrm{Def}_G(V,R)$ of \emph{deformations} of $V$ over $R$ 
to be the set of isomorphism classes of lifts of $V$ over $R$.
The \emph{deformation functor} $F_V:\cat\to \mathrm{Sets}$ is defined to be the 
covariant functor which sends a ring $R$ in $\cat$ to $\mathrm{Def}_G(V,R)$ and
a morphism $\alpha:R\to R'$  in $\hat{\mathcal{C}}$ to the set map
\begin{eqnarray*}
F_V(\alpha):\qquad \mathrm{Def}_G(V,R) & \to & \mathrm{Def}_G(V,R')\\
  \; [M,\phi] &\mapsto&  [R'\otimes_{R,\alpha}M,\phi_\alpha]
\end{eqnarray*}
where $\phi_\alpha$ is the composition 
$k\otimes_{R'}(R'\otimes_{R,\alpha}M)\cong k\otimes_RM\xrightarrow{\phi} V$.

Sometimes it is useful to describe the deformation functor $F_V$ in terms of matrix groups.
Choosing a $k$-basis of $V$, we can identify $V$ with $k^n$ where $n=\mathrm{dim}_k\,V$.
The $G$-action on $V$ is then given by
a continuous homomorphism $\overline{\rho}:G \to \mathrm{GL}_n(k)$, which 
is called a \emph{residual representation}.
Let $R$ be a ring in $\hat{\mathcal{C}}$ and denote the reduction map
$\mathrm{GL}_n(R)\to \mathrm{GL}_n(k)$ induced by $\pi_R:R\to k$ also by $\pi_R$.
By a \emph{lift} of $\overline{\rho}$ over $R$ we mean a continuous
homomorphism $\rho:G\to \mathrm{GL}_n(R)$ such that 
$\pi_R\circ\rho=\overline{\rho}$.
Such a lift defines a $G$-action on $M=R^n$, and with the obvious
isomorphism $\phi:k\otimes_RM\to V$ such a lift defines a deformation
$[M,\phi]$ of $V$ over $R$.
Two lifts $\rho,\rho':G\to \mathrm{GL}_n(R)$ of $\overline{\rho}$ over $R$ give rise
to the same deformation if and only if they are \emph{strictly equivalent},
that is, if one can be brought into the other by conjugation by a matrix in the kernel of $\pi_R$.
In this way, the choice of a basis of $V$ gives rise to an identification of
$\mathrm{Def}_G(V,R)$ with the set $\Def_G(\overline{\rho},R)$ of strict 
equivalence classes of lifts of $\overline{\rho}$ over $R$. In particular, the 
deformation functor $F_{\overline{\rho}}:\cat\to\mathrm{Sets}$ associated to
strict equivalence classes of lifts of $\overline{\rho}$ is naturally isomorphic to the
deformation functor $F_V$. Note that if $\alpha:R\to R'$ is a morphism in $\cat$, then  
$F_{\overline{\rho}}(\alpha)$ is the set map $\mathrm{Def}_G(\overline{\rho},R) 
\to \mathrm{Def}_G(\overline{\rho},R')$ which sends $[\rho]$ to $[\alpha\circ\rho]$, 
where we denote the morphism $\mathrm{GL}_n(R)\to \mathrm{GL}_n(R')$ induced by 
$\alpha$ also by $\alpha$. In the following, we identify $F_V=F_{\overline{\rho}}$.

\subsection{Universal deformation rings}
\label{ss:udr}

Assume the notation from the previous subsection.
The functor $F_V$ is representable if there is a ring $R(G,V)$
in $\hat{\mathcal{C}}$ and a lift $(U(G,V),\phi_U)$ of $V$ over $R(G,V)$
such that for all $R$ in $\hat{\mathcal{C}}$ the map
\begin{eqnarray*}
f_R:\qquad \mathrm{Hom}_{\hat{\mathcal{C}}}(R(G,V),R) &\to& 
\mathrm{Def}_G(V,R)\\
\;\alpha &\mapsto& F_V(\alpha)([U(G,V),\phi_U])
\end{eqnarray*}
is bijective. Put differently, this is the case if and only if $F_V$ is naturally isomorphic to the
Hom functor $\mathrm{Hom}_{\cat}(R(G,V),-)$.
If this is the case then $R(G,V)$ is called the \emph{universal
deformation ring} of $V$ and $[U(G,V),\phi_U]$ is called the \emph{universal
deformation} of $V$. The defining property determines $R(G,V)$ 
and $[U(G,V),\phi_U]$ uniquely up to a unique isomorphism.

A slightly weaker notion can be useful if the functor $F_V$ is
not representable.  The ring $k[\epsilon]$ of dual numbers with $\epsilon^2 = 0$
has a $W$-algebra structure such that the maximal ideal $\mathfrak{m}_W$ of $W$ annihilates
$k[\epsilon]$.  One says $R(G,V)$ is a \emph{versal deformation ring} of $V$ if
the maps $f_R$ are surjective for all $R$, and bijective
for $R=k[\epsilon]$.  These conditions determine $R(G,V)$ uniquely up to
isomorphism, but the isomorphism need not be unique.

\begin{theorem} {\rm (\cite[\S 1.2]{maz1},  \cite[Prop. 7.1]{desmitlenstra})}
\label{thm:mazurudr}
Suppose $G$ satisfies the following $p$-finiteness condition:
\\[.2cm]
$\mathbf{(\Phi_p)}$ 
For every open subgroup $J$ of finite index in $G$, there
are only finitely many continuous homomorphisms from $J$ to $\mathbb{Z}/p$.
\\[.2cm]
Then every finite dimensional continuous representation $V$ of 
$G$ over $k$ has a versal deformation ring.  If $\mathrm{End}_{kG}(V)=k$, then $V$ has a 
universal deformation ring.
\end{theorem}

To prove the existence of versal deformation rings, 
Mazur verified Schlessinger's criteria \cite{Sch} of pro-representability
of Artin functors for the deformation functor $F_V$. He also proved that $F_V$ 
is continuous, meaning that for all objects $R$ in $\hat{\mathcal{C}}$
we have $\displaystyle F_V(R)=\lim_{\longleftarrow} F_V(R/\mathfrak{m}_R^i)$.
Note that  Mazur assumed $k$ to be a finite field in \cite{maz1}. 
However, his proofs in \cite[\S 1.2]{maz1} go through in the more general case we are considering.

In the case when $\mathrm{End}_{kG}(V)=k$,
de Smit and Lenstra took a different approach in \cite{desmitlenstra}
which proceeds in three main steps: First they let $G$ be finite and considered the functor which
assigns to each $R$ a certain set of homomorphisms $G\to\mathrm{GL}_n(R)$.
They showed that this functor is representable by defining the corresponding universal ring
by generators and relations. Taking projective limits, they obtained a similar result for
profinite $G$. Finally, they concluded the construction by passing to a suitable closed subring,
which is either generated by the traces of the elements of $G$ if $V$ is absolutely irreducible,
or by a larger collection of elements as suggested by Faltings if $\mathrm{End}_{kG}(V)=k$.

In the case when $V$ is absolutely irreducible, another approach was given by
Rouquier in \cite{rouquier}, using pseudo-characters and the results in 
\cite{nyssen}. The main point in this latter construction is that being a pseudo-character function
has a universal solution, providing another description of the universal deformation ring in
this case.

\subsection{Pseudocompact modules}
\label{ss:pseudocompact}

In this subsection, we briefly describe another viewpoint of lifts and deformations using pseudocompact
modules. This is useful, for example, when generalizing deformations of group representations to
deformations of objects in derived categories (see for example \cite{bcderived1,bcderived2}).
Pseudocompact rings, algebras and modules have been studied, for example, in 
\cite{ga1,ga2,brumer}, which also serve as references for the following statements.

Assume the notation from subsection \ref{ss:deffunctor}.
For $R \in \mathrm{Ob}(\hat{\mathcal{C}})$, let $R[[G]]$ be the completed group algebra of the usual
abstract group algebra $RG$ of $G$ over  $R$, i.e. $R[[G]]$ is the projective limit of the usual group 
algebras $R[G/U]$ as $U$ ranges over the open normal subgroups of $G$ (where we put brackets
around $G/U$ for better readability). Giving a finitely generated free $R$-module $M$ on which $G$
acts continuously is the same as giving a topological $R[[G]]$-module $M$ 
which is finitely generated and free as an $R$-module.

The completed group algebra $R[[G]]$ is a so-called \emph{pseudocompact ring}, i.e. it is
a complete Hausdorff topological ring which admits a basis of open neighborhoods of $0$ consisting of
two-sided ideals $J$ for which $R[[G]]/J$ is an Artinian ring. 
In particular, $R[[G]]$ is the projective limit of Artinian quotient rings having the discrete topology.
Since $R$ is a commutative pseudocompact ring and $R[[G]]$ is an $R$-algebra
and since the open neighborhood basis of $0$ can be chosen to consist of two-sided ideals $J$ 
for which $R[[G]]/J$ has finite length as $R$-module, $R[[G]]$ is moreover a \emph{pseudocompact 
$R$-algebra}.
A complete Hausdorff topological $R[[G]]$-module $M$ is said to be a \emph{pseudocompact 
$R[[G]]$-module}
if $M$ has a basis of open neighborhoods  of $0$ consisting of submodules $N$ for which
$M/N$ has finite length. 
It follows that an $R[[G]]$-module is pseudocompact if and only if it is the projective limit of 
$R[[G]]$-modules of finite length having the discrete topology. 
We denote the category of pseudocompact $R[[G]]$-modules by $\mathrm{PCMod}(R[[G]])$.
Note that $\mathrm{PCMod}(R[[G]])$ is an abelian category with exact projective limits. 

A pseudocompact $R[[G]]$-module $M$ is said to be \emph{topologically free} on a set
$X=\{x_i\}_{i\in I}$ if $M$ is isomorphic to the product of a family $(R[[G]]_i)_{i\in I}$ where
$R[[G]]_i=R[[G]]$ for all $i$. In particular, a topologically free pseudocompact $R[[G]]$-module
on a finite set is the same as a finitely generated abstractly free $R[[G]]$-module.

As before, assume $V$ is a finite dimensional vector space over $k$
with a continuous $k$-linear action of $G$ on $V$, which is given by 
a continuous group homomorphism from $G$ to the discrete group $\Aut_k(V)$.
Then $V$ is  a pseudocompact $k[[G]]$-module. Moreover, any lift of $V$
over a ring $R$ in $\hat{\mathcal{C}}$ is given by a pseudocompact $R[[G]]$-module $M$
which is finitely generated and abstractly free as an $R$-module together with an isomorphism 
$\phi: k \otimes_R M \to V$ in $\mathrm{PCMod}(k[[G]])$. Note that in principle, we should use the
completed tensor product  $\hat{\otimes}_R-$ (see \cite[\S2]{brumer}) rather than the usual tensor product
$\otimes_R$. However, since $k$ is finitely generated as a pseudocompact $R$-module, it follows 
that the functors $k\otimes_R -$ and $k\hat{\otimes}_R-$ are naturally isomorphic.

Every topologically free pseudocompact $R[[G]]$-module is a projective object in 
$\mathrm{PCMod}(R[[G]])$, and every pseudocompact
$R[[G]]$-module is the quotient of a topologically free $R[[G]]$-module. Hence
$\mathrm{PCMod}(R[[G]])$ has enough projective objects. 
If $M$ and $N$ are 
pseudocompact $R[[G]]$-modules, then we define the right derived functors
$\mathrm{Ext}^n_{R[[G]]}(M,N)$ by using a projective resolution of $M$. 

For all positive integers $n$, we have an isomorphism of $k$-vector spaces
\begin{equation}
\label{eq:cohomology}
\HH^n(G,\mathrm{Hom}_k(V,V))\cong \mathrm{Ext}^n_{k[[G]]}(V,V)
\end{equation}
where $\HH^n$ is ``continuous'' cohomology.
Since $V$ is finite dimensional and discrete, this isomorphism can be deduced by looking at direct limits
over all open normal subgroups $U$ of $G$ which act trivially on $V$.

\subsection{Tangent space}
\label{ss:tangent}

We continue to assume the notation from subsection \ref{ss:deffunctor}. Moreover, we
assume that $G$ satisfies the $p$-finiteness condition $\mathbf{(\Phi_p)}$ from
Theorem \ref{thm:mazurudr}.
For simplicity, we also assume that $V$ has a universal deformation ring $R(G,V)$. 
In other words,
the deformation functor $F_V$ is represented by $R(G,V)$ and thus naturally isomorphic to 
the Hom functor $\mathrm{Hom}_{\cat}(R(G,V),-)$. 

The \emph{tangent space} of $F_V$ is defined as 
$$t_{F_V}=F_V(k[\epsilon])\cong \mathrm{Hom}_{\cat}(R(G,V),k[\epsilon])$$
where, as before,
$k[\epsilon]$ is the ring of dual numbers with $\epsilon^2 = 0$. Since in $\cat$ there is
only one morphism $R(G,V)\to k$, we obtain an isomorphism
$$\mathrm{Hom}_{\cat}(R(G,V),k[\epsilon])\xrightarrow{\cong} \mathrm{Hom}_k\left(
\frac{\mathfrak{m}_{R(G,V)}}{\mathfrak{m}_{R(G,V)}^2+pR(G,V)}\,,\,k\right)$$
of $k$-vector spaces. For $R$ in $\cat$, we call 
$$t_R^*=\frac{\mathfrak{m}_R}{\mathfrak{m}_R^2+pR}$$
the Zariski cotangent space of $R$. Hence we obtain an isomorphism of tangent spaces
$$t_{F_V}\cong \mathrm{Hom}_k(t_{R(G,V)}^*,k)=t_{R(G,V)}.$$
The following result gives a connection of $t_{F_V}$ to the cohomology of 
$\mathrm{Hom}_k(V,V)$.

\begin{proposition} {\rm (\cite[\S1.2]{maz1}, \cite[\S 22]{maz2})}
\label{prop:tangentspace}
There is a natural isomorphism of $k$-vector spaces
$$t_{F_V}\cong \HH^1(G,\mathrm{Hom}_k(V,V))\cong \mathrm{Ext}^1_{k[[G]]}(V,V).$$
If $r=\mathrm{dim}_k\mathrm{Ext}^1_{k[[G]]}(V,V)$ then $R(G,V)$ is isomorphic to 
a quotient algebra of the power series algebra $W[[t_1,\ldots,t_r]]$ in $r$ commuting
variables and $r$ is minimal with this property.
\end{proposition}

As in the previous subsection, $\mathrm{Ext}^1_{k[[G]]}(V,V)$ means $\mathrm{Ext}^1$ in the 
category of pseudocompact $k[[G]]$-modules; see also Equation (\ref{eq:cohomology}).
The main idea of the proof is to notice that if $(M,\phi)$ is a lift of $V$ over $k[\epsilon]$
then by restricting the scalars from $k[\epsilon]$ to $k$ we may view $M$ as a 
$k$-vector space of dimension $2\cdot\mathrm{dim}_kV$, with a $k$-linear continuous action
of $G$. Identifying the $k[[G]]$-modules $\epsilon M$ and $M/\epsilon M$ with $V$ (using
$\phi$), we then see $M$ as an extension of $V$ by $V$ in the category of pseudocompact 
$k[[G]]$-modules:
$$\mathcal{E}:\quad 0\to V \xrightarrow{\iota} M \xrightarrow{\tau} V\to 0.$$
Sending the element of $t_{F_V}$ corresponding to the isomorphism class 
$[M,\phi]$ to the element of $\mathrm{Ext}^1_{k[[G]]}(V,V)$ corresponding to $\mathcal{E}$,
we obtain a well-defined map
$$s:t_{F_V}\to \mathrm{Ext}^1_{k[[G]]}(V,V)$$
which is a $k$-vector space homomorphism. The inverse map of $s$ is obtained by
going backward:
Given an extension $\mathcal{E}$ as above, we define a $k[\epsilon]$-structure on $M$
by letting $\epsilon$ act as the composition $\iota\circ \tau$, enabling us to view
$M$ as a lift of $V$ over $k[\epsilon]$.

\subsection{Obstructions}
\label{ss:obstruct}

We continue to assume the notation from subsection \ref{ss:deffunctor},
that $G$ satisfies $\mathbf{(\Phi_p)}$ and that $V$ has a 
universal deformation ring $R(G,V)$. Let $\overline{\rho}:G\to\mathrm{GL}_n(k)$
be a residual representation corresponding to $V$.

Suppose we have a surjective morphism $R_1\to R_0$ in $\cat$
and assume that the kernel $I$ satisfies $I\cdot \mathfrak{m}_{R_1}=0$,
so $I$ has the structure of a $k$-vector space. Suppose we have a lift
$\rho_0:G\to\mathrm{GL}_n(R_0)$ of $\overline{\rho}$ over $R_0$. We want to describe
the obstruction to lifting $\rho_0$ to $R_1$. 

Let $\gamma_1:G\to\mathrm{GL}_n(R_1)$ be a set-theoretic lift of $\rho_0$. Since
$\gamma_1$ is a homomorphism modulo $I$, we obtain a  2-cocycle
\begin{eqnarray*}
c:\quad G\times G &\to& 1+\mathrm{Mat}_n(I)\\
(g_1,g_2) &\mapsto& \gamma_1(g_1g_2)\gamma_1(g_2)^{-1}\gamma_1(g_1)^{-1}.
\end{eqnarray*}
Identifying the multiplicative group $1+\mathrm{Mat}_n(I)\subset \mathrm{GL}_n(R_1)$ 
with the additive group
$\mathrm{Mat}_n(I)^+$, we obtain the following isomorphisms:
$$1+\mathrm{Mat}_n(I)\cong \mathrm{Mat}_n(I)^+\cong\mathrm{Mat}_n(k)\otimes_kI
\cong \mathrm{Hom}_k(V,V)\otimes_kI.$$
Here the action of $G$ on $\mathrm{Mat}_n(k)$, respectively $\mathrm{Hom}_k(V,V)$,
is given by the usual conjugation action; note that this is also called
the \emph{adjoint representation} of $\overline{\rho}$ and denoted by
$\mathrm{ad}(\overline{\rho})$.
If we replace $\gamma_1$ by a different
set-theoretic lift, this changes $c$ by a 2-coboundary. Therefore we obtain an
element
$$[c]\in\HH^2(G,\mathrm{Hom}_k(V,V)\otimes_kI)\cong 
\HH^2(G,\mathrm{Hom}_k(V,V))\otimes_kI$$
which gives the \emph{obstruction} to lifting $\rho_0$ to $R_1$. The class
$[c]$ is sometimes called the \emph{obstruction class} of $\rho_0$ relative to
the morphism $R_1\to R_0$. 

\begin{proposition} {\rm (\cite[\S1.6]{maz1}, \cite[Thm. 2.4]{bockle})}
\label{prop:obstructions}
If $r=\mathrm{dim}_k\mathrm{Ext}^1_{k[[G]]}(V,V)$ and 
$s=\mathrm{dim}_k\mathrm{Ext}^2_{k[[G]]}(V,V)$ then $R(G,V)$ is isomorphic to 
a quotient algebra $W[[t_1,\ldots,t_r]]/J$ where $s$ is an upper bound on the minimal
number of generators of $J$.
\end{proposition}

The main idea of the proof is as follows: Let $\rho_u:G\to\mathrm{GL}_n(R(G,V))$ be
a universal lift of $\overline{\rho}$.
Let $\hat{R}=W[[t_1,\ldots,t_r]]$ and consider the sequence
$$0\to J'=J/(J\mathfrak{m}_{\hat{R}})\to R'=\hat{R}/(J\mathfrak{m}_{\hat{R}}) \to R(G,V)\to 0.$$
Since $\mathfrak{m}_{R'}$ annihilates $J'$, we have an obstruction
$[c]\in \HH^2(G,\mathrm{Hom}_k(V,V))\otimes_kJ'$ to lifting $\rho_u$ to $R'$. This obstruction
depends only on the strict equivalence class of $\rho_u$ and not on the chosen representation
$\rho_u$.  One then shows that the $k$-linear map
\begin{eqnarray*}
\mathrm{Hom}_k(J',k) &\to& \HH^2(G,\mathrm{Hom}_k(V,V)) \cong \mathrm{Ext}^2_{k[[G]]}(V,V)\\
f &\mapsto& (1\otimes f)([c])
\end{eqnarray*}
is injective, which implies the proposition.

\subsection{Deformation rings in number theory}
\label{ss:explicitudr}

One of the main uses of deformation rings in number theory has been to establish
a relationship between certain Galois representations and automorphic forms.
More precisely, let $\overline{\rho}:G_\mathbb{Q}\to\mathrm{GL}_2(k)$ be an
absolutely irreducible representation, where $k$ is a finite field of characteristic $p>2$ and 
$G_\mathbb{Q}$ is the absolute Galois group of $\mathbb{Q}$. 
Suppose that $\overline{\rho}$ is modular in the sense that it corresponds to a modular form 
(modulo $p$) which is an eigenfunction of Hecke operators. The idea is to prove that
all reasonable lifts of $\overline{\rho}$ to $p$-adic representations are modular by
establishing an isomorphism between a universal deformation ring, which parameterizes
lifts of $\overline{\rho}$ with bounded ramification and satisfying appropriate deformation conditions, 
and a Hecke algebra, which parameterizes certain
lifts of $\overline{\rho}$ which are modular of some fixed level. 

Taylor and Wiles established such an isomorphism in \cite{wita}, which then led to the
proof of Fermat's Last Theorem in \cite{wi}. The Taylor-Wiles method has been further
refined by many people, such as Diamond \cite{diamond} and Fujiwara \cite{fujiwara}.
In \cite{khare2003}, Khare gave an alternative approach 
for semistable $\overline{\rho}$ by establishing an isomorphism 
$R_Q\cong T_Q$, where $Q$ is a so-called auxiliary set of primes, $R_Q$ is the 
universal deformation ring for 
lifts of $\overline{\rho}$ minimally ramified away from $Q$ and satisfying appropriate deformation
conditions and $T_Q$ is the analogous 
Hecke algebra.

Apart from the proof of Fermat's Last Theorem in \cite{wita,wi} and the 
proof of the general Shimura-Taniyama-Weil conjecture 
in \cite{bcdt}, deformation rings also played
an important role in the proof of Serre's modularity conjecture \cite{serremod} 
by Khare-Wintenberger \cite{khawin} and Kisin \cite{kisin2},
which asserts that every absolutely irreducible representation 
$\overline{\rho}:G_\mathbb{Q}\to\mathrm{GL}_2(k)$ with odd determinant
is modular.

\medskip

Suppose $k$ is a finite field of characteristic $p>2$, $K$ is a number field,
$S$ is a finite set of primes of $K$ containing the primes over $p$ and the infinite ones, and
$\overline{\rho}:G_K\to\mathrm{GL}_2(k)$ is an absolutely irreducible representation 
unramified outside $S$. 
It is often desirable to have an explicit presentation of a universal
deformation ring $R$, which parameterizes lifts of $\overline{\rho}$ satisfying
certain deformation conditions, in terms of a power series algebra over $W=W(k)$ modulo an ideal
given by a (minimal) number of generators. In the following, we describe 
several results of B\"ockle in this respect, some of which also played an important role
in \cite{khawin}.

Since the relations occurring in the universal deformation ring $R$
often come from the obstructions of the associated local deformation problems 
$\overline{\rho}_{\mathfrak{p}}:G_{K_{\mathfrak{p}}}\to\mathrm{GL}_2(k)$ for $\mathfrak{p}\in S$,
B\"ockle considered in \cite{bodemuskin}  the problem of finding the universal deformation, or a 
smooth cover of it, in the local case where the relevant pro-$p$ group is an arbitrary 
Demu\v{s}kin group. He showed that the corresponding universal deformation ring
is a complete intersection, flat over $W$, and with the 
(minimal) number of generators given by the $k$-dimension of 
$\HH^2(G_{K_{\mathfrak{p}}},\mathrm{ad}_{\overline{\rho}})$. 
Moreover, he applied his local results to the global situation.
For example, he
gave conditions under which the universal deformation ring of an odd, absolutely irreducible representation $G_{\mathbb{Q}}\to \mathrm{GL}_2(k)$, unramified outside
$S$, can be described explicitly,
thus generalizing a result of Boston \cite{boston}.

In \cite{bockle,bopresent}, B\"ockle studied in more detail the connection between local and global 
deformation functors.  In \cite{bockle},
he presented a rather general class of (global) deformation functors of $\overline{\rho}$ 
that satisfy local deformation conditions and investigated for those, under what conditions the 
global deformation functor is determined by the local deformation functors
corresponding to primes $\mathfrak{p}\in S$. 
B\"ockle gave precise conditions under which the local functors are 
sufficient to describe the global functor. 
These conditions involve the vanishing of a second Shafarevich-Tate group
and auxiliary primes as introduced by Taylor and Wiles in \cite{wita}.
In \cite{bopresent}, B\"ockle provided generalizations and simplified proofs for some of the 
results in \cite{bockle}.


\section{Universal deformation rings of modules for finite groups}
\label{s:udrfinite}

Assume the notation from subsection \ref{ss:deffunctor}.
If $\mathrm{End}_{kG}(V)=k$, the construction of the universal deformation ring
$R(G,V)$ by de Smit and Lenstra in \cite{desmitlenstra} shows that 
$R(G,V)$ is the inverse limit of the universal deformation rings $R(H,V)$ when $H$ ranges over 
all  finite discrete quotients of $G$ through which the $G$-action on $V$ factors. Thus to answer 
questions about the ring structure of $R(G,V)$, it is natural to first consider the case when $G=H$ 
is finite. 

For the remainder of this paper, we assume that $G$ is finite. The representation theory of $kG$
when $p$ divides $\# G$ is very beautiful but difficult. To avoid rationality questions and to simplify notation,
we assume that $k$ is algebraically closed. More precisely, we make the following assumptions.

\begin{hypothesis}
\label{hypo:finite}
Let $k$ be an algebraically closed field of positive characteristic $p>0$, let $W=W(k)$
be the ring of infinite Witt vectors over $k$, let $G$ be a \textbf{finite} group, and
let $V$ be a finitely generated $kG$-module.
\end{hypothesis}

It follows as before that $V$ has a universal deformation ring if its endomorphism
ring $\mathrm{End}_{kG}(V)$ is isomorphic to $k$. Note that
$\mathrm{End}_{kG}(V)\cong \HH^0(G,\mathrm{Hom}_k(V,V))$.
When $G$ is finite, Tate cohomology
groups often play an important role. Therefore, the question arises if we can
use the 0-th Tate cohomology group,
$$\hat{\HH}^0(G,\mathrm{Hom}_k(V,V)) = \mathrm{End}_{kG}(V)/(s_G\cdot \mathrm{End}_k(V))$$
where $s_G=\sum_{g\in G} g$, to obtain a criterion for the existence of a universal deformation
ring of $V$. 
Let $\mathrm{PEnd}_{kG}(V)$ denote the ideal of $\mathrm{End}_{kG}(V)$
consisting of all $kG$-module endomorphisms of $V$ factoring through a projective $kG$-module.
Then $\mathrm{PEnd}_{kG}(V)$ is equal to $s_G\cdot \mathrm{End}_k(V)$, which implies that 
$\hat{\HH}^0(G,\mathrm{Hom}_k(V,V)) = \mathrm{End}_{kG}(V)/\mathrm{PEnd}_{kG}(V)$.
The quotient ring $\mathrm{End}_{kG}(V)/\mathrm{PEnd}_{kG}(V)$
is called the \emph{stable endomorphism ring of $V$} and is
denoted by $\underline{\mathrm{End}}_{kG}(V)$.
Note that in general $\mathrm{PEnd}_{kG}(V)\neq 0$, i.e. $\mathrm{End}_{kG}(V)$ properly
surjects onto $\underline{\mathrm{End}}_{kG}(V)$.
We have the following result:

\begin{proposition} {\rm (\cite[Prop. 2.1]{bc}, \cite[Rem. 2.1]{3quat})}
\label{prop:blchin} 
Assume Hypothesis $\ref{hypo:finite}$ and that
the stable endomorphism ring
$\underline{\mathrm{End}}_{kG}(V)$ is isomorphic to $k$. Then $V$ has a 
universal deformation ring. Moreover, if $R$ is in $\hat{\mathcal{C}}$ and
$(M,\phi)$ and $(M',\phi')$ are lifts of $~V$ over $R$ such that $M$ and $M'$ are
isomorphic as $RG$-modules then $[M,\phi]=[M',\phi']$.
\end{proposition}

The main point of the proof of this proposition is to show that if
$(M,\phi)$ is a lift of $V$ over an Artinian object $R$ in $\cat$, then the ring homomorphism 
$R \to\underline{\mathrm{End}}_{RG}(M)$ coming from the action of $R$ on $M$ via scalar 
multiplication is  surjective. 

\begin{remark}
\label{rem:nophi}
The last statement of Proposition \ref{prop:blchin} means that the particular $kG$-module
isomorphism $k\otimes_RM\xrightarrow{\phi} V$ from a lift $(M,\phi)$ is not important
when  $\underline{\mathrm{End}}_{kG}(V)\cong k$,
which significantly simplifies computations. 

Note that this is not true in general, as can be seen in the following example.
Let $G=\langle \sigma\rangle$ be a cyclic group of order $p$ and let $V=k\oplus k$
with trivial $G$-action.
Let $M$ be the $k[\epsilon]G$-module with $M=k[\epsilon]\oplus k[\epsilon]$ and
$\sigma$ acting as $\left[\begin{array}{cc}1&\epsilon\\ 0&1\end{array}\right]$. 
Consider $\phi,\phi':k\otimes_{k[\epsilon]}M\to V$ with
$\phi=\left[\begin{array}{cc}1&0\\ 0&1\end{array}\right]$ and
$\phi'=\left[\begin{array}{cc}0&1\\1&0\end{array}\right]$. Then $(M,\phi)$ and
$(M,\phi')$ are two non-isomorphic lifts of $V$ over $k[\epsilon]$.
\end{remark}

The following result analyzes $R(G,V)$ further in the case when 
$\underline{\mathrm{End}}_{kG}(V)\cong k$. Here $\Omega$ denotes the syzygy functor
or Heller operator,
i.e. if $\pi:P_V\to V$ is a projective $k$G-module cover of $V$ then $\Omega(V)$ denotes the 
kernel of $\pi$ (see, for example, \cite[\S20]{alp}).

\begin{proposition} 
\label{prop:defhelp}
{\rm (\cite[Cors. 2.5 and 2.8]{bc}).}
Assume Hypothesis $\ref{hypo:finite}$ and that the stable endomorphism ring
$\underline{\mathrm{End}}_{kG}(V)$ is isomorphic to $k$.
\begin{enumerate}
\item[(i)] Then $\underline{\mathrm{End}}_{kG}(\Omega(V))\cong k$, and $R(G,V)$ and 
$R(G,\Omega(V))$ are isomorphic.
\item[(ii)] There is a non-projective indecomposable $kG$-module $V_0$ $($unique up to isomorphism$)$ such that $\underline{\mathrm{End}}_{kG}(V_0)\cong k$, $V$ is isomorphic to $V_0\oplus Q$ for some projective $kG$-module $Q$, and $R(G,V)$ and $R(G,V_0)$ are isomorphic.
\end{enumerate}
\end{proposition}

The main ideas of the proof are as follows: Since the deformation functor $F_V$ is continuous,
most of the arguments can be carried out for the restriction of $F_V$ to the full subcategory
$\mathcal{C}$ of $\hat{\mathcal{C}}$ of Artinian objects. For part (i), one shows that 
the syzygy functor $\Omega$ induces an isomorphism between the restrictions of the functors
$F_V$ and $F_{\Omega(V)}$ to   $\mathcal{C}$. For part (ii), one uses that the projective
$kG$-module $Q$ can be lifted to a projective $RG$-module $Q_R$ for
every $R$ in $\mathcal{C}$  to show that
there is an isomorphism between the restrictions of the functors
$F_V$ and $F_{V_0}$ to $\mathcal{C}$.

\medskip

Recall that $kG$ can be written as a finite direct product of blocks
$$kG=B_1\times \cdots \times B_r$$
where the blocks $B_1,\ldots,B_r$ are in one-to-one correspondence with the primitive central
idempotents of $kG$. (For a good introduction to block theory, we refer the reader to \cite[Chap. IV]{alp}.)
If $B$ is a block of $kG$, there is associated to it a conjugacy class
of $p$-subgroups of $G$, called the defect groups of $B$. The defect groups
measure how far $B$ is away from being a full matrix ring; they also determine the
representation type of $B$. More precisely, $B$ has finite representation type if and only
if its defect groups are cyclic; $B$ has infinite tame representation type if and only if
$p=2$ and the defect groups of $B$ are dihedral, semi-dihedral or generalized quaternion;
and $B$ has wild representation type in all other cases. (This result follows from \cite{hi,br,bd}. 
A description of this result together with an introduction to the representation type can also be 
found in \cite[Intro. and Sect. I.4]{erd}.)

Proposition \ref{prop:defhelp}(ii) says that if the stable endomorphism ring of $V$ is isomorphic
to $k$ and we want to determine the universal deformation ring $R(G,V)$ then we may
assume that $V$ is non-projective indecomposable. But then $V$ belongs to a unique block
of $kG$, and we can use the theory of blocks, as introduced by Brauer and developed by
many other authors, to determine the universal deformation ring $R(G,V)$. 


\section{Brauer's generalized decomposition numbers and universal deformation rings}
\label{s:brauer}

We continue to assume Hypothesis \ref{hypo:finite}. Our goal in this section is to show how
Brauer's generalized decomposition numbers can be used in certain cases to determine
the isomorphism type of the universal deformation ring $R(G,V)$. We first give
a brief introduction to these generalized decomposition numbers.

\subsection{Brauer's generalized decomposition numbers}
\label{ss:Brauergendec}

The usual decomposition numbers were introduced by Brauer and Nesbitt in 
\cite{brauernesbitt1937} (see also \cite{brauernesbitt}). They allow us to express
the values of the ordinary irreducible characters of $G$ on $p$-regular 
elements of $G$, i.e. elements of order prime to $p$, by means  of the absolutely irreducible
$p$-modular characters of $G$. More precisely, if $\zeta_1, \zeta_2,\ldots$
are the ordinary irreducible characters of $G$ and $\varphi_1,\varphi_2,\ldots$ are the
absolutely irreducible $p$-modular characters of $G$, then we have a formula
\begin{equation}
\label{eq:usualdec}
\zeta_\mu(g) = \sum_\nu d_{\mu\nu}\varphi_\nu(g)
\end{equation}
provided $g$ is a $p$-regular element of $G$. 
The $d_{\mu\nu}$ are non-negative integers, called the 
\emph{decomposition numbers of $G$ for $p$}.
As Brauer wrote in \cite[p. 192]{brauerordinaryandmodular}:

\smallbreak
{``We may say that the group characters $\zeta_\mu$ of $G$ are built up by the modular 
characters $\varphi_\nu$, and it is possible to obtain a deeper insight into the nature of the 
ordinary group characters by the use of the modular characters and their properties. However, 
it is disturbing that we have to restrict ourselves to $p$-regular elements.''}
\smallbreak

For this reason, Brauer introduced generalized decomposition numbers in 
\cite{brauerordinaryandmodular}. The value $\zeta_\mu(g)$ on an element 
$g\in G$ whose order is divisible by $p$ is then expressed by means of the absolutely irreducible
$p$-modular characters of certain subgroups $C_i$ of $G$. The corresponding 
generalized decomposition numbers $d^i_{\mu\nu}$ are not necessarily rational
integers, but they
are algebraic integers in a cyclotomic field of $p$-power order roots of unity.
More precisely, Brauer defined $d^i_{\mu\nu}$ as follows.

Suppose $\# G=p^a m'$ where $m'$ is relatively prime to $p$, and let $P$ be a fixed 
Sylow $p$-subgroup of $G$.
Let $u_0 = 1, u_1, u_2, \ldots,u_h$ be a complete system of
representatives of $G$-conjugacy classes of $p$-power order elements in $G$ with $u_i\in P$
for all $1\le i\le h$. Every  conjugacy class of $G$ contains an element of the form
$u_i v$ where $i\in\{0,1,\ldots, h\}$ is uniquely determined by the class and $v$ is a $p$-regular element in the
centralizer $C_G(u_i)$.
For each $0\le i\le h$, let $v_{i,1},\ldots,v_{i,\ell_i}$ be a complete system of representatives
of $C_G(u_i)$-conjugacy classes of $p$-regular elements in $C_G(u_i)$ with $v_{i,1}=1$.
Then $\{u_iv_{i,j}\;|\; 0\le i\le h,1\le j\le \ell_i\}$ is a complete set of representatives
of the conjugacy classes of $G$. 
Moreover, for each $0\le i\le h$, there are precisely
$\ell_i$ absolutely irreducible $p$-modular characters of $C_G(u_i)$, which we denote by
$\varphi_1^i,\ldots,\varphi_{\ell_i}^i$. As before, let $\zeta_1,\zeta_2,\ldots$ be the
ordinary irreducible characters of $G$. Then 
\begin{equation}
\label{eq:brauer1}
\zeta_\mu(u_iv_{i,j})= \sum_{\nu=1}^{\ell_i} d_{\mu\nu}^i \,\varphi_\nu^i(v_{i,j})
\end{equation}
for all $0\le i\le h$, $1\le j\le \ell_i$. 
The $d_{\mu\nu}^i$ are called the \emph{generalized decomposition numbers of $G$}. 
For $i = 0$, we have $u_0=1$ and $C_G(u_0)=G$, and the $d_{\mu\nu}^i$
coincide with the usual decomposition numbers $d_{\mu\nu}$ of $G$ in Equation (\ref{eq:usualdec}).
In general, $d_{\mu\nu}^i$ is an algebraic integer in the field of the $p^{\alpha_i}$-th roots of
unity where $p^{\alpha_i}$ is the order of $u_i$.
In particular, $d_{\mu\nu}^i$ can be viewed to belong to $W[\omega_i]$ if $\omega_i$ is a 
primitive $p^{\alpha_i}$-th root of unity.
In \cite[Sect. 6]{brauerdarst2}, Brauer moreover showed that if $\zeta_\mu$ belongs
to the block $B$ of $kG$, then the generalized decomposition number
$d_{\mu\nu}^i$ vanishes if $\varphi_\nu^i$ belongs to a block of $kC_G(u_i)$ whose Brauer 
correspondent in $G$ is not equal to $B$.

\subsection{Universal deformation rings of certain modules belonging to infinite tame blocks}
\label{ss:udrtame}

We now focus on a certain class of modules belonging to blocks of infinite tame representation
type for which Brauer's generalized decomposition numbers can be used to determine
their universal deformation rings. This subsection is based on the paper \cite{brauerpaper},
and details can be found there.
Recall from section \ref{s:udrfinite} that a block $B$ of $kG$
has infinite tame representation type if and only if $p=2$ and the defect groups of $B$ are
either dihedral, or semi-dihedral, or generalized quaternion.

In \cite{brIV,brauer2,olsson}, Brauer and Olsson determined the generalized decomposition
numbers for all the ordinary irreducible characters belonging to infinite tame blocks.
Moreover, they proved that an infinite tame block has at most three isomorphism classes of
simple modules.
In \cite{erd}, Erdmann classified all infinite tame blocks up to Morita equivalence
by providing a list of quivers and relations for their basic algebras. 

We make the following assumptions.

\begin{hypothesis}
\label{hyp:alltheway}
Assume Hypothesis $\ref{hypo:finite}$. Additionally, assume that $p=2$, $V$ is
indecomposable with $\underline{\mathrm{End}}_{kG}(V)\cong k$, and that 
$V$ belongs to a non-local block $B$ of $kG$ of infinite tame representation type
with a defect group $D$ of order $2^n$. 
Let $F$ be the fraction field of $W$, and let 
$\overline{F}$ be
a fixed algebraic closure of $F$.
\end{hypothesis}

We want to concentrate on those $V$ for which Brauer's generalized decomposition
numbers carry the most information. 
More precisely, we call a module $V$ as in Hypothesis \ref{hyp:alltheway}
\emph{maximally ordinary} if the $2$-modular character of $V$
is the restriction to the $2$-regular conjugacy classes of an ordinary irreducible character $\chi$
such that for every $\sigma\in D$ of maximal $2$-power order,
Brauer's generalized decomposition numbers corresponding to $\sigma$ and $\chi$
do not all lie in $\{0,\pm 1\}$. In other words, using the notation of Equation (\ref{eq:brauer1}),
if $\chi=\zeta_\mu$ and $\sigma$ is conjugate
in $G$ to $u_i$, then there exists an absolutely
irreducible $2$-modular character $\varphi_\nu^i$ of $C_G(u_i)$
such that $d_{\mu\nu}^i \not\in\{0,\pm 1\}$.

By \cite{brauer2,olsson}, 
there are precisely $2^{n-2}-1$ ordinary irreducible characters of height 1 belonging to 
$B$ if $n\ge 4$. Moreover, they all define the same $2$-modular character when they are
restricted to the $2$-regular conjugacy classes.
If $n=3$, then there are either 1 or 3 ordinary irreducible characters of 
height 1 belonging to $B$, depending on whether $D$ is dihedral or quaternion.
If $n=2$, then there are no ordinary irreducible characters of height 1 belonging to $B$.
Recall that the height of an ordinary irreducible character $\chi$ belonging to $B$
is $b-a+n$, where $2^a$ (resp. $2^b$) is the maximal $2$-power dividing
$\# G$ (resp. $\mathrm{deg}(\chi)$). Since $n$ is the defect of the block $B$, it
follows that $b-a+n$ is a non-negative integer (see, for example, \cite[Sect. 56.E and Cor. (57.19)]{CR}).

Suppose $n\ge 4$. By \cite{brauer2,olsson}, exactly one of the $2^{n-2}-1$ ordinary irreducible characters of 
height 1 belonging to $B$ is realizable over $F$, i.e. it
corresponds to an absolutely irreducible $FG$-module.
Moreover, the remaining $2^{n-2}-2$ characters of height 1 
are precisely the ordinary irreducible characters belonging to $B$ for which
the generalized decomposition numbers
corresponding to maximal $2$-power order elements in $D$ do not all lie in $\{0,\pm 1\}$.
We have the following result:

\begin{theorem}
{\rm (\cite[Thm. 1.1 and Cor. 6.2]{brauerpaper}).}
\label{thm:supermain}
Assume Hypothesis $\ref{hyp:alltheway}$.
Then $V$ is maximally ordinary if and only if $n\ge 4$
and the $2$-modular character of $V$ is equal to the restriction
to the $2$-regular conjugacy classes
of an ordinary irreducible character of $G$ of height $1$. 
Suppose $V$ is maximally ordinary.
There exists a monic polynomial $q_n(t)\in W[t]$ of degree $2^{n-2}-1$
which depends on $D$ but not on $V$ and which can be given explicitly
such that  either
\begin{enumerate}
\item[(i)] $R(G,V)/2R(G,V)\cong k[[t]]/(t^{2^{n-2}-1})$, in which case $R(G,V)$ is isomorphic to
$W[[t]]/(q_n(t))$, or
\item[(ii)] $R(G,V)/2R(G,V)\cong k[[t]]/(t^{2^{n-2}})$, in which case  $R(G,V)$ is isomorphic to
$W[[t]]/(t\,q_n(t),2\,q_n(t))$.
\end{enumerate}
In all cases, the ring $R(G,V)$ is isomorphic to a subquotient ring of $WD$, and it is
a complete intersection if and only if we are in case $(i)$.
\end{theorem}

A precise description of the maximally ordinary modules $V$ belonging to $B$ is given in 
\cite[Lemma 6.1 and Cor. 6.2]{brauerpaper}.
A formula for the polynomials $q_n(t)$ can be found
in \cite[Def. 5.3 and Rem. 5.4]{brauerpaper}.

We now discuss the main ideas of the proof of Theorem \ref{thm:supermain}.
For details we refer the reader to \cite{brauerpaper}.
The first statement of the theorem follows from the results in \cite{brIV,brauer2,olsson}.
Suppose now that $n\ge 4$.
As noted above,  there are then precisely $2^{n-2}-1$ ordinary irreducible 
characters of height 1 belonging to $B$. Moreover, these characters fall into $n-2$ 
Galois orbits under the action of 
$\mathrm{Gal}(\overline{F}/F)$:
$$\mathcal{F}_0,\mathcal{F}_1,\ldots,\mathcal{F}_{n-3}$$
where $\#\mathcal{F}_j=2^j$ for $0\le j\le n-3$. If $\xi$ is the ordinary character
which is the sum of all the characters of height 1, then $\xi$ can be realized by
an $FG$-module
$$X=X_0\oplus X_1\oplus \cdots \oplus X_{n-3}$$
where each $X_j$ is a simple $FG$-module with Schur index 1 corresponding
to the orbit $\mathcal{O}_j$.

The main steps to prove Theorem \ref{thm:supermain} are  as follows:
Suppose $V$ is maximally ordinary.
First, we use Erdmann's description of the basic algebra of $B$ to show that
$R(G,V)$ is a quotient algebra of $W[[t]]$ and that $V$ can be lifted to
$k[[t]]/(t^{n-2}-1)$. Moreover, we show that this lift is given by an indecomposable
$B$-module $\overline{U'}$ of $V$ such that
$\mathrm{End}_{kG}(\overline{U'})\cong k[[t]]/(t^{n-2}-1)$. 
Next, we use the usual decomposition numbers, together with the description of the
projective indecomposable $B$-modules and \cite[Prop. (23.7)]{CR}
and \cite[Lemma 2.3.2]{3sim} to show that $\overline{U'}$ can be lifted to $W$.
Moreover, we show that this lift is given by an indecomposable $WG$-module
$U'$ which is free over $W$ with $F\otimes_WU'\cong X$.
Then we use Brauer's generalized decomposition numbers to show that
$\mathrm{End}_{WG}(U')\cong W[[t]]/(q_n(t))$ and that $U'$ is free as a module
for $\mathrm{End}_{WG}(U')$. This then implies that $W[[t]]/(q_n(t))$
is a quotient ring of the universal deformation ring $R(G,V)$.

To complete the proof of Theorem \ref{thm:supermain}, we use
again Erdmann's description of the basic algebra of $B$ to determine the 
universal mod 2 deformation ring $R(G,V)/2 R(G,V)$. It follows that 
the isomorphism type of $R(G,V)/2 R(G,V)$ depends on whether or not the
stable Auslander-Reiten quiver $\Gamma_s(B)$ of $B$ contains 3-tubes. 
Note that if $D$ is dihedral then $\Gamma_s(B)$ always contains 3-tubes, 
whereas if $D$ is generalized quaternion then $\Gamma_s(B)$ never contains 3-tubes, 
and if $D$ is semi-dihedral then $\Gamma_s(B)$ may or may not contain 3-tubes. 
We show that $R(G,V)/2R(G,V)\cong k[[t]]/(t^{n-2}-1)$ if $\Gamma_s(B)$
contains no 3-tubes, and that $R(G,V)/2R(G,V)\cong k[[t]]/(t^{n-2})$ if
$\Gamma_s(B)$ does contain 3-tubes. 
In the first case, the universal deformation ring of $V$ is 
$R(G,V)\cong W[[t]]/(q_n(t))$, whereas in the second case we use
\cite[Lemma 2.3.3]{3sim} to show that
$R(G,V)\cong W[[t]]/(t\,q_n(t),2\,q_n(t))$.

\bibliographystyle{amsplain}

\end{document}